\documentclass[reqno,a4paper]{amsart}
\usepackage[active]{srcltx}
\usepackage{fancyhdr}
\usepackage[utf8]{inputenc}
\DeclareRobustCommand{\SkipTocEntry}[5]{} 
\usepackage{times}
\usepackage{latexsym,amsopn,amssymb,amsmath,mathabx}
\changenotsign
\usepackage{amsthm}
\usepackage{amsfonts,amsbsy,amscd,stmaryrd,bbold}
\usepackage{dsfont}
\usepackage{paralist}
\usepackage{hyperref}
\hypersetup{colorlinks=true,linkcolor=Emerald,urlcolor=blue} 

\usepackage[usenames,dvipsnames]{color}

\usepackage{mathrsfs} 

\newcommand{\C}{\mathbb{C}}
\newcommand{\N}{\mathbb{N}}
\newcommand{\R}{\mathbb{R}}
\newcommand{\Z}{\mathbb{Z}}

\newcommand{\cd}{\mathcal{D}}
\newcommand{\ce}{\mathcal{E}}

\newcommand{\cg}{\mathcal{G}}
\newcommand{\ch}{\mathcal{H}}

\newcommand{\ck}{\mathcal{K}}
\newcommand{\cl}{\mathcal{L}}

\newcommand{\cq}{\mathcal{Q}}
\newcommand{\cR}{\mathcal{R}}
\newcommand{\cs}{\mathcal{S}}
\newcommand{\ct}{\mathcal{T}}
\newcommand{\cu}{\mathcal{U}}

\def\rmi{\mathrm{i}}

\def\jap#1{\langle {#1} \rangle}

\DeclareMathOperator*{\w*lim}{w*-lim}

\def\qed{\hfill $\Box$\medskip}

\newtheorem{theorem}{Theorem}[section]

\newtheorem{proposition}[theorem]{Proposition}

\newtheorem{assumption}[theorem]{Assumption}

\newtheorem{definition}[theorem]{Definition}

\numberwithin{equation}{section}

\begin{document}

\title[ ]{On the Limiting Absorption Principle for Schr\"odinger operators on waveguides
  \\[2mm]
 {\tiny \today}
}

\author[A. Martin]{Alexandre Martin} \address{A. Martin,
  D\'epartement de Math\'ematiques, Universit\'e de Cergy-Pontoise,
  95000 Cergy-Pontoise, France}
\email{alexandre.martin@u-cergy.fr }

\begin{abstract}
We prove a Limiting Absorption Principle for Schr\"odinger operators in tubes about infinite curves embedded in the Euclidian space with different types of boundary conditions. The argument is based on the Mourre theory with conjugate operators different from the generator of dilations which is usually used in this case, and permits to prove a Limiting Absorption Principle for Schr\"odinger operators in singular waveguides.
\end{abstract}

\maketitle
\tableofcontents

\section{Introduction}

The purpose of this article is to prove a limiting absorption principle for a certain class of Schr\"odinger operators on a waveguide and to study the nature of their essential spectrum. To do this, we will use a general technique due to E. Mourre \cite{Mo1} which involves a local version of the positive commutator method due to C.R. Putnam \cite{Pu1,Pu2}. If we want to use this theory to prove a limiting absorption principle for a self-adjoint operator $H$, the idea is to search for a second operator $A$, which is taken in general self-adjoint, such that $H$ is regular with respect to $A$ in a certain sense and such that $H$ satisfies the Mourre estimate on a set $I$ in the following sense
\[E(I)[H,iA]E(I)\geq c_0E(I)+K\]
where $E(I)$ is the spectral measure of $H$ on $I$, $c_0>0$ and $K$ is a compact operator.

Let $n\in\N$, $n\geq2$ and $\Sigma\subset\R^{n-1}$ an open bounded set. Consider the straight tube $\Omega=\R\times\Sigma$. In this article, we will study many types of boundary conditions, mainly the Dirichlet and Neumann conditions. Thus we will suppose that the boundary of $\Sigma$ is of class $C^1$, but in certain cases, this assumption is not necessary (for Dirichlet boundary conditions for example). When we want to apply the Mourre theory to Schr\"odinger operators on waveguides, we usually used the following operator
\[A=A_D^y\otimes\mathbb{1}_\Sigma=\frac{y\partial_y+\partial_y y}{2i},\]
with domain $C^\infty_c(\Omega)$.

This operator is a generator of dilations only in the unbounded direction of the waveguide. We can easily show that this operator is essentially self-adjoint, that, if $\Delta$ is a Laplacian on the waveguide (with Dirichlet, Neumann, Robin boundary conditions), $\Delta\in C^\infty(A)$ and that, if $V$ is the operator of multiplication by a function of class $C^1$ then $[V,iA]=-q_y\partial_{y}V$ (see \cite{Ben,DDI1,DDI2,KT}). 

Using this conjugate operator for the Dirichlet Laplacian, we can see the following

\begin{theorem}[Theorem 2.16 de \cite{KT}]\label{th: KT}
Let $\Sigma$ an open bounded connected set of $\R^{n-1}$, $n\geq2$, and denote by $\ct$ the set of eigenvalues of $\Delta_D^\Sigma$, the Dirichlet Laplacian on $\Sigma$. Let $\Omega=\R\times\Sigma$ and $H=\Delta_D+V$ on $L^2(\Omega)$ with Dirichlet conditions on the boundary and with $V$ the operator of multiplication by a real function. Assume that $V$ satisfy the following:
\begin{itemize}
\item $V\in L^\infty(\Omega)$;

\item $\lim\limits_{R\rightarrow\infty} \sup\limits_{x\in(\R\backslash[-R,R])\times\Sigma}|V(x)|=0$;

\item There is $\theta>0$ and $C>0$ such that $|\partial_y V(x)|\leq C (1+y^2)^{-\frac{1+\theta}{2}},\forall x=(y,\sigma)\in\Omega$.  
\end{itemize}
Alors 
\begin{enumerate}[(i)]
\item $\sigma_{ess}(H)=[\kappa,\infty)$ with $\kappa=\inf \ct$;

\item $\sigma_{sc}(H)=\emptyset$;

\item $\sigma_p(H)\cup\ct$ is closed and countable.

\item $\sigma_p(H)\backslash \ct$ is composed of finitely degenerated eigenvalues, which can accumulate at the points of $\ct$ only;

\item The limit $R(\lambda\pm \rmi 0)=\w*lim\limits_{\mu\rightarrow0} R(\lambda\pm\rmi\mu)$ exists, locally uniformly in $\lambda\in (\kappa,\infty)\backslash\ct$ outside of the eigenvalues of $H$, with $R(z)=(H-z)^{-1}$.
\end{enumerate}
\end{theorem}
A similar result can be proved with Neumann boudary conditions by taking $\ct$ the set of eigenvalues of the Neumann Laplacian on $\Sigma$. In this case, $\kappa=0$.

In this article, we will study different types of conjugate operators to prove a Limiting Absorption principle. In section \ref{s: A_D toute direction}, we will see why a generator of dilations in all directions does not seem to be a good choice of conjugate operator to use the Mourre theory. In section \ref{s: A_u guide d'onde}, we will see how the choice of a conjugate operator with decay in the momentum variable permits to use the Mourre theory for Schr\"odinger operators on curved waveguide with singular potential when we want to prove a Limiting Absorption Principle far from thresholds. We will also see how to prove a Limiting Absorption Principle near thresholds.

\section{A generator of dilations in all directions}\label{s: A_D toute direction}

As it was said in the introduction, the use of a generator of dilations $A$ only in the unbounded direction of the waveguide cause a problem near the eigenvalues of the Laplacian on $\Sigma$. We can think that this problem is due to the fact that, in the commutator between the Laplacian and $A$, the Laplacian does not appear in its entirely. An idea to solve that may be to take as conjugate operator a generator of dilations in all directions, like it is the case for Schr\"odinger operators on the Euclidian space. Here we will see that this choice of conjugate operator seems to be a bad choice.

\subsection{Results}
Let $\Sigma$ be an open bounded connected set of $\R^{n-1}$ such that $0\in\Sigma$, with a sufficiently regular boundary (we can suppose that the boundary is of class $C^1$ for example). Let $\Omega=\R\times\Sigma$ be a waveguide. We consider the operator $A_0=-i 2^{-1}(x\cdot \nabla+\nabla\cdot x)$ the generator of dilations in all directions with domain $\cd(A_0)=C^\infty_c(\Omega)$. Now we will give some properties of this operator and its relation with Dirichlet, Neumann and Robin Laplacians.

\begin{proposition}\label{prop: A_D guide d'onde }
\begin{enumerate}
\item $A_0$ is a symmetric operator without self-adjoint extensions;

\smallskip

\item Neumann and Robin Laplacians are not of class $C^1(\bar{A_0})$;

\smallskip

\item Dirichlet Laplacian is of class $C^1(\bar{A_0})$ but not of class $C^1_u(\bar{A_0})$. Moreover, if we denote $(\lambda_k)_{k\in\N^*}$ eigenvalues of the Dirichlet Laplacian on $\Sigma$, numbered in ascending order, then for all $k\in\N^*$, if $I\subset]\lambda_k,\lambda_{k+1}[$ and $|I|\leq \inf\limits_{1\leq j\leq k}\lambda_{j+1}-\lambda_j$, then the Mourre estimate is true on $I$ with $\bar{A_0}$ as conjugate operator.
\end{enumerate}
\end{proposition}

Thus we can not use the Mourre theory for Neumann and Robin Laplacian. We can moreover remark that, in certain case, the Mourre estimate is satisfied for the Dirichlet Laplacian for intervals $I$ with size larger than $\inf\limits_{1\leq k\leq n}\lambda_{k+1}-\lambda_k$. For example, if $\Sigma=[-1,1]$, this assumption on the size of $I$ can be replaced by $|I|\leq \inf\limits_{1\leq k\leq n-1}\lambda_{k+2}-\lambda_k$.

Now we will show the different results.

To begin, we will prove that $A_0$ is symmetric. Let $f,g\in\cd(A_0)$. we have:
 \begin{eqnarray*}
 (f,A_0g)&=&\left(f,\frac{x\cdot \nabla+\nabla\cdot x}{2i}g\right)\\
 &=&-\left(\frac{1}{2i}f,(x\cdot \nabla+\nabla\cdot x)g\right).
 \end{eqnarray*}
 To simplify notations, let $F=\frac{1}{2i}f$.
 \begin{eqnarray}\label{eq:symetrique}
 (f,A_0g)&=&-(F,(x\cdot \nabla+\nabla\cdot x)g)\nonumber\\
 &=&-\int_\Omega\bar{F}(x)\left(x\cdot \nabla g(x)+\nabla\cdot xg(x)\right)dx\nonumber\\
 &=&\int_\Omega\left(x\cdot \nabla \bar{F}(x)+\nabla\cdot x\bar{F}(x)\right)g(x)dx+\sum_{k=1}^n\int_{\partial\Omega}x_k \bar{F}(x)g(x)\nonumber\\
 &=& (A_0 f,g).
 \end{eqnarray}
 Then $A_0$ is symmetric. Moreover, we can remark that the assumption $f\in\cd(A_0)$ is not necessary. In fact, under the assumption $x\cdot\nabla f\in L^2(\Omega)$, the previous computation stays true, boundary terms disappearing by  Dirichlet conditions which are assumed for $g$. Thus, we have $\cd(A_0^*)\supset\{u\in L^2(\Omega), x\cdot \nabla u\in L^2(\Omega)\}$. Moreover, by the definition of the adjoint, we can remark that if $u\in\cd(A_0^*)$, then $u$ is necesseraly in the set $L^2(\Omega)$ and, by \eqref{eq:symetrique}, we can show the equality $\cd(A_0^*)=\{u\in L^2(\Omega), x\cdot \nabla u\in L^2(\Omega)\}$. We can remark that if we have other conditions on the boundary of the waveguide (Neumann, Robin,...) for $A_0$, this operator can not be symmetric because of the addition of boundary terms in \eqref{eq:symetrique}. 
 
 Now we will see if $A_0$ has a self-adjoint extension. For this, we will use the caracterisation of symmetric operators given in the Corollary at the begining of page 141 of \cite{RS2}:
\begin{proposition}\label{prop: max sym RS}
Let $A$ a symmetric closed operator.
Denote $n_\pm=dim Ker(A^*\mp \rmi Id)$ its deficiency index.
Then
\begin{enumerate}
\item $A$ is self-adjoint if and only if $n_+=n_-=0$;

\item $A$ has a self-adjoint extension if and only if  $n_+=n_-$;

\item If $n_+=0\not=n_-$ or $n_-=0\not=n_+$, then $A$ has no nontrivial symmetric extensions ($A$ is called maximal symmetric).
\end{enumerate}
\end{proposition}
Moreover, by Theorem X.1 of \cite{RS2}, we know that, for $A$ a closed symmetric operator, the dimension of spaces $Ker(A^*-\lambda ID)$ is constant throughout the open upper half-plane $\{z\in\C,\Im(z)>0\}$ and throughout the open lower half-plane $\{z\in\C,\Im(z)<0\}$. Since $A_0$ is symmetric and densly defined, we know that $A_0^*=\bar{A_0}^*$ (see \cite[Theorem VIII.1]{RS1}). Thus we have to search for the dimension of $Ker(A_0^*-\rmi\lambda)$ for different values of $\lambda$ ($\lambda>0$ or $\lambda<0$).

We begin by the case $\lambda>0$. For $\phi\in \cd(A_0)$, let $\tilde{\phi}$ its extension by $0$ to $\R^n$. $\tilde{\phi}\in C^\infty_c(\R^n)$.
For $t>0$, denote $\phi_t(x)=\tilde{\phi}(tx)$. Since $\Sigma$ is an open set which contains $0$, for $\phi\in \cd(A_0)$, there exist $t_0\geq1$ such that  $\phi_t\in \cd(A_0)$ if  $t\geq t_0$.

We know that $\phi_t$ converge simply almost everywhere to $\psi\equiv0$ when $t$ goes to $+\infty$. In particular, for all $f\in L^2(\Omega)$, $(f,\phi_t)$ goes to $0$ when $t$ goes to $+\infty$ by dominated convergence. Let $\lambda>0$ and $f\in L^2(\Omega)$ such that 
\[\left(\frac{x\cdot \nabla+\nabla\cdot x}{2i}-i\lambda\right)f=0.\]
In particular, $x\cdot \nabla f=-(\frac{n}{2}+\lambda)f\in L^2(\Omega)$. Moreover $f\in\cd(A_0^*)$. Thus, for $\phi\in C^\infty_c(\Omega)$,  we have for $t\geq t_0$
\begin{eqnarray*}
\partial_t(f,\phi_t)&=&(f,\partial_t \phi_t)\\
&=&(f,x\cdot \nabla \phi)(tx))\\
&=& \frac{1}{t}\left(f,x\cdot\nabla(\phi_t)(x)\right)\\
&=& \frac{1}{t}\left(f,(\rmi A_0-\frac{n}{2})\phi_t\right)\\
&=& -\frac{1}{t}\left((\rmi A_0^*+\frac{n}{2})f,\phi_t\right)\\
&=&\frac{1}{t}(\lambda-\frac{n}{2})(f,\phi_t).
\end{eqnarray*}
Therefore, for $t\geq t_0$,
\[(f,\phi_t)=(f,\phi)\exp\left((\lambda-\frac{n}{2})\ln(t)\right).\]
by taking the limit when $t$ goes to $+\infty$, we deduce that for $\lambda\geq\frac{n}{2}$, $(f,\phi)=0$ for all $\phi\in C^\infty_c(\Omega)$ and then $f=0$.

Thus, we have shown that $Ker(A_0^*-i\lambda I)=\{0\}$ for $\lambda\geq\frac{n}{2}$. Therefore, $Ker(\bar{A_0}^*-i\lambda I)=\{0\}$ for $\lambda\geq\frac{n}{2}$. By \cite[Theorem X.1]{RS2} applied to the closure of $A_0$,we deduce that $Ker(\bar{A_0}^*-i\lambda I)=\{0\}$ for all $\lambda>0$.

Now, we look to the case $\lambda<0$. Let $f$ be a solution of the equation $A_0^* f=i\lambda f$. We can write this equation with the form
\[x\cdot\nabla f=-(\lambda+\frac{n}{2})f.\]
By composing with the unitary operator of convertion into polar coordinates, denoting $(\theta_i)_{i=1,\cdots,n-1}$ angular variables, we can see that this equation can be written $r\partial_r f=-(\lambda+\frac{n}{2})f$. Thus, for $\lambda<-1-\frac{n}{2}$ and $C:]-\pi,\pi]^{n-1}\rightarrow\R$ a function $C^\infty$ with compact support included in  $]-\pi,\pi[^{n-1}\backslash\{0\}$, this equation admits for solution the function

\[h_\lambda(r,\theta_1,\cdots,\theta_{n-1})=C(\theta_1,\cdots,\theta_{n-1})r^{-\lambda-\frac{n}{2}}.\]

Remark that $r\partial_r h_\lambda(r,\theta_1,\cdots,\theta_{n-1})=-(\lambda+\frac{n}{2}) h_\lambda(r,\theta_1,\cdots,\theta_{n-1})$. Moreover, since $\lambda<-1-\frac{n}{2}$ and $C$ is bounded, $h_\lambda\in C^1(\R^n)$. Since $C$ has a compact support in \\$]-\pi,\pi[^{n-1}\backslash\{0\}$, we can also remark that $h_\lambda\in L^2(\Omega)$ which implies that $x\cdot\nabla h_\lambda=r\partial_r h_\lambda\in L^2(\Omega)$. Thus, we have $h_\lambda\in \cd(A_0^*)$. Therefore $h_\lambda\in Ker(A_0^*-i\lambda)$ which implies that $dim Ker(\bar{A_0}^*+i)>0$. Thus, we are in the setting of the point (3) of Proposition \ref{prop: max sym RS} which implies that $\bar{A_0}$ is maximal symmetric.

The point (1) of Proposition \ref{prop: A_D guide d'onde } is thus proved.

With $\bar{A_0}$ as conjugate operator, we can not use the classic Mourre theory but a Mourre theory adapted to maximal symmetric operators (see \cite{GGM}). 

Now, we will see what happen for Laplacian with different type of boundary conditions. We will denote $\Delta_D$ the Dirichlet Laplacian, $\Delta_N$ the Neumann Laplacian and if computations are similar for all types of Laplacian, we denote it $\Delta$. To apply this Mourre theory, adapted to maximal symmetric operators, to a Laplacian (with Dirichlet, Neumann or Robin conditions on the boundary) with $\bar{A_0}$ as conjugate operator, it is necessary to have $\Delta \in C^1(A_0)$. Now we recall a caracterisation of this regularity adapted to our context

\begin{proposition}[Proposition 2.22 of \cite{GGM}]\label{prop: 2.22 de GGM}
Let $S$ a self adjoint operator on $\ch$ and $A$ a maximal symmetric operator on $\ch$. Then $H\in C^1(A)$ if and only if the two following conditions are satisfied:
\begin{enumerate}
\item There exists $c\geq0$ such that for all $u\in\cd(A^*)\cap\cd(S)$ and $v\in\cd(A)\cap\cd(S)$, $|(u,[S,A]v)|\leq c\|u\|_S\|v\|_S$,

\smallskip

\item There exist $z\in\rho(S)$ such that $\{f\in\cd(A), R(z)f\in\cd(A)\}$ is a core for $A$ and $\{f\in\cd(A^*), R(\bar{z})f\in\cd(A^*)\}$ is a core for $A^*$. 
\end{enumerate}
\end{proposition}
To simplify notations, let $A_1=\bar{A_0}$ the closure of $A_0$. 
Let $z\in\rho(\Delta )$.

Let $u\in \cd(A_1^*)$. Let $v=(\Delta -\bar{z})^{-1}u$. Thus $v$ satisfied $(\Delta-\bar{z})v=u$ with Dirichlet, Neumann or Robin conditions on the boundary according to the Laplacian considerated. To prove that $v\in\cd(A_1^*)$, it suffices to show that $x\cdot\nabla v\in L^2$. By definition of $v$, we have
\begin{eqnarray}\label{eq: prop 2.3 dernier point}
x\cdot\nabla v&=&x\cdot\nabla (\Delta-\bar{z})^{-1}u\nonumber\\
&=&x\cdot (\Delta-\bar{z})^{-1}\nabla u\nonumber\\
&=& (\Delta-\bar{z})^{-1}x\cdot\nabla u+[x,(\Delta-\bar{z})^{-1}]\cdot\nabla u\nonumber\\
&=&(\Delta-\bar{z})^{-1}x\cdot\nabla u-(\Delta-\bar{z})^{-1}[x,\Delta]\cdot (\Delta-\bar{z})^{-1}\nabla u\nonumber\\
&=&(\Delta-\bar{z})^{-1}x\cdot\nabla u-2i(\Delta-\bar{z})^{-2} \Delta u.
\end{eqnarray}
By assumptions, $x\cdot\nabla u\in L^2$. Thus $v\in\cd(A_1^*)$. In particular, $\{f\in\cd(A_1^*), (\Delta-\bar{z})^{-1}f\in\cd(A_1^*)\}=\cd(A_1^*)$ and is a core for $A_1^*$.

\subsection{The case of Neumann and Robin Laplacian}
Now, we will give the proof of point (2) of Proposition \ref{prop: A_D guide d'onde }. To show that, we will try to apply Proposition \ref{prop: 2.22 de GGM}. Remark that, by \eqref{eq: prop 2.3 dernier point}, the second part of assumption (2) of Proposition \ref{prop: 2.22 de GGM} is already proved. Thus, we will show that other conditions are not satisfied. We begin with the Neumann Laplacian.

Let $z\in\rho(\Delta_N)$ and $\ce=\{f\in\cd(A_1), (\Delta_N-z)^{-1}f\in\cd(A_1)\}$.
Let $u\in\ce$. Let $v\in L^2$ such that $v=(\Delta_N-z)^{-1}u$. Then $v$ satisfied
\[\begin{cases}
\Delta_N v-zv=u\quad \text{ dans }\R\times\Sigma\\
\frac{\partial v}{\partial n}=0 \quad \text{ sur }\R\times\partial\Sigma
\end{cases}\]
where $\frac{\partial}{\partial n}$ is the normal derivative.
Moreover $v\in\cd(A_1)$. Thus $v\restriction_{\R\times\partial\Sigma}=0$.
For $f\in L^2(\R\times\Sigma)$, denote $\hat{f}$ the Fourier transform of $f$ with respect to the first variable $y$. Therefore $v$ satisfies:
\[\begin{cases}
\Delta_\Sigma \hat{v}+(\xi^2-z)\hat{v}=\hat{u}\quad \text{ dans }\R\times\Sigma\\
\frac{\partial \hat{v}}{\partial n}=0 \quad \text{ sur }\R\times\partial\Sigma
\end{cases}.\]
Let $\alpha(\xi)\in\C$ such that $(\alpha(\xi))^2=\xi^2-\bar{z}$ and let
\[w(\xi,\sigma)=\exp\left(\frac{\alpha(\xi)}{(n-1)^{1/2}}\sum_{k=1}^{n-1}\sigma_k\right).\]
We can remark that $w$ satisfies $\Delta_\Sigma w+(\xi^2-\bar{z})w=0$. Moreover, since $\Sigma$ is bounded, for all $\xi\in\R$, $w(\xi,\cdot)\in L^2(\Sigma)$. Then, we can define the linear map $L:L^2(\Omega)\rightarrow L^1(\R)$ by $L(f)(\xi)=\int_\Sigma \hat{f}(\xi,\sigma)w(\xi,\sigma)d\sigma$. By Green formula, we have for all $u\in\ce$, $L(u)=0$. This implies that $\ce\subset L^{-1}(\{0\})$. But, if $g\in C^\infty_c(\Sigma)$, $g\not=0$ and $g\geq0$, by denoting $u_1(y,\sigma)=\exp(-\frac{y^2}{2})g(\sigma)$, then $\hat{u}_1(\xi,\sigma)=u_1(\xi,\sigma)\geq0$ for all $(\xi,\sigma)\in\Omega$. Thus we obtain $L(u_1)>0$ and moreover $u_1\in\cd(A_1)$. Since $L$ is continuous, we deduce that $\ce$ can not be dense into $\cd(A_1)$. Thus it is not a core for $A_1$ which already proved that $\Delta_N\notin C^1(A_1)$. 

Remark that if we replace Neumann boundary conditions by Robin boundary conditions, the same result appears. In fact, since $v\in\cd(A_1)$, $v$ satisfies Dirichlet boundary conditions. Asking that $v$ satisfies Neumann boundary conditions is then equivalent to the fact of asking that $v$ satisfies Robin boundary conditions. Thus the Robin Laplacian is not of class $C^1(A_1)$.

For Neumann Laplacian, we can still search to know if the commutator is bounded from $\cd(\Delta_N)$ to its dual space. To do this, we will make the computation of the commutator between $\Delta_N$ and $A_1$. To simplify computations, we will see only the case $\Sigma=\prod_{k=1}^{n-1}[a_i,b_i]$.

Let $f\in\cd(\Delta_N)\cap\cd(A_1)$. We have:

 \begin{eqnarray*}
(f,[\Delta_N,iA_1]f)&=&(f,-2\partial_y^2 f)+\sum_{k=1}^{n-1}\int_{\R\times\Sigma}(-\partial_{\sigma_k}^2 \bar{f})(\sigma_k\partial_{\sigma_k} f+\frac{f}{2})d\sigma\\
& &+\int_{\R\times\Sigma}(-\partial_{\sigma_k}^2 f)(\sigma_k\partial_{\sigma_k} \bar{f}+\frac{\bar{f}}{2})d\sigma.
\end{eqnarray*}
We can remark that since $\Sigma$ is a rectangle, the outwardly normal vector of $\Sigma$ is a vector of the standrad basis (or the oposite of a vector from the basis). By Fubini Theorem and by integration by part, we obtain 
 \[(f,[\Delta_N,iA_1]f)=2\left((f,-\partial_y^2 f)+\sum_{k=1}^{n-1}(f,-\partial_{\sigma_k}^2 f)\right)=-2\int_{\R\times\Sigma}\bar{\nabla f}\nabla f dyd\sigma.\]
 
 The form domain of the commutator is equal to $\cd(\Delta_N)\cap\cd(A_1)$. But, since $C^\infty_c\subset\cd(\Delta_N)\subset\ch^1$ and $C^\infty_c\subset\cd(A_1)\subset\{f\in L^2, f\restriction_{\R\times\partial\Sigma}=0\}$, we can deduce that
 \[C^\infty_c\subset\cd(\Delta_N)\cap\cd(A_1)\subset\ch^1_0.\]
  Thus we have
 \[(f,[\Delta_N,iA_1]f)=2(f,\Delta_Df).\]
 Therefore, the first commutator is positive and one would think that the strict Mourre estimate is true on all intervall. Unfortunately, the Dirichlet Laplacian is not bounded with respect to the Neumann Laplacian ($\cd(\Delta_N)\nsubset\cd(\Delta_D)$), and then condition (1) of Proposition \ref{prop: 2.22 de GGM} is not satisfied. 
 
Since the commutator is positive, we can ask if we can use Lemma 3.12 of \cite{GGM} for which it is not necessary for $H$ to be regular with respect to $A$ but only regular with respect to an operator $H'$ such that there is $c>0$ with
 \[H'+c\jap{H}\geq\jap{H},\]
  where $H'=[\Delta_N,iA_1]$. In fact, since $H'=2\Delta_D$, $[\Delta_N,iH']=0$ with appropriate domain, we can show that the condition (1) of Proposition \ref{prop: 2.22 de GGM} is satisfied with $A=H'$. On the other hand, as previously, we can show that $\{f\in\cd(\Delta_D),(\Delta_N+z)^{-1}f\in\cd(\Delta_D)\}$ is not a core for $\Delta_D$. Thus, the condition (2) of Proposition \ref{prop: 2.22 de GGM} is not satisfied. Therefore $\Delta_N\notin C^1(H')$ wich prevent us to use Lemma 3.12 from \cite{GGM}. 

\subsection{The case of Dirichlet Laplacian}
Now, we will give the proof of point (3) of Proposition \ref{prop: A_D guide d'onde }. To show that, we will try to apply Proposition \ref{prop: 2.22 de GGM}. Remark that, by \eqref{eq: prop 2.3 dernier point}, the second part of assumption (2) of Proposition \ref{prop: 2.22 de GGM} is already proved. Thus, we will show that other conditions are satisfied. 

Let $(y,\sigma)\in\R\times\Sigma$ a point of the waveguide, $A_D^y$ the generator of dilations in the direction $y$ and $A_D^\sigma$ the generator of dilations in the direction $\sigma$ with Dirichlet boundary conditions. We can remark that $A_1$ can be written $A_1=A_D^y\otimes \mathbb{1}_\Sigma +\mathbb{1}_\R\otimes A_D^\sigma$. We can remark also that, since $\Sigma$ is bounded, $\cd(A_D^\sigma)=\ch^1_0(\Sigma)$. Let $u\in\cd(A_D^y\otimes \mathbb{1}_\Sigma)\cap\cd(\mathbb{1}_\R\otimes A_D^\sigma)$ and $w=(\Delta_D-z)^{-1}u$. Since $A_D^y\otimes \mathbb{1}_\Sigma$ is self-adjoint, as in \eqref{eq: prop 2.3 dernier point}, we can show that $w\in \cd(A_D^y\otimes \mathbb{1}_\Sigma)$. Moreover, $w\in\cd(\Delta_D)\subset L^2(\R,\ch^1_0(\Sigma))$. Thus $w\in\cd(\mathbb{1}_\R\otimes A_D^\sigma)$. Therefore, since $\cd(A_D^y\otimes \mathbb{1}_\Sigma)\cap\cd(\mathbb{1}_\R\otimes A_D^\sigma)\subset \cd(A_1)$, we have
\[C^\infty_c\subset \cd(A_D^y\otimes \mathbb{1}_\Sigma)\cap\cd(\mathbb{1}_\R\otimes A_D^\sigma)\subset\{f\in\cd(A_1), (\Delta_D-z)^{-1}f\in\cd(A_1)\}.\]
Thus $\{f\in\cd(A_1), (\Delta_D-z)^{-1}f\in\cd(A_1)\}$ is a core for $A_1$.

To show that $\Delta_D\in C^1(A_1)$ using Proposition \ref{prop: 2.22 de GGM}, it remains to show that the commutator is bounded from the domain of $\Delta_D$ into its dual space.

To begin, we can remark
\[\Delta_D=-\partial_y^2\otimes\mathbb{1}_\Sigma+\mathbb{1}_\R\otimes(\Delta_D^\Sigma),\]
where $\Delta_D^\Sigma$ is the Laplacian on $\Sigma$ with Dirichlet boundary conditions.
In particular, since our conjugate operator has the same form, we have, in the sense of sesquilinear form on $\cd(\Delta_D)\cap\cd(A_1)$:
\begin{eqnarray*}
[\Delta_D,iA_1]&=&[-\partial_y^2,iA_D^y]\otimes\mathbb{1}_\Sigma+\mathbb{1}_\R\otimes[\Delta_D^\Sigma,iA_D^\sigma]\\
&=&-2\partial_y^2\otimes\mathbb{1}_\Sigma+\mathbb{1}_\R\otimes[\Delta_D^\Sigma,iA_D^\sigma].
\end{eqnarray*}
Thus, it remains to compute $[\Delta_D^\Sigma,iA_D^\sigma]$.
Let $f\in \cd(\Delta_D^\Sigma)\cap\cd(A_D^\sigma)$. We have:
\begin{eqnarray*}
(f,[\Delta_D^\Sigma,iA_D^\sigma]f)&=&(\Delta_D^\Sigma f,iA_D^\sigma f)+(iA_D^\sigma f,\Delta_D^\Sigma f)\\
&=& \int_\Sigma(\Delta_\Sigma \bar{f})\frac{\sigma\cdot\nabla_\sigma f+\nabla_\sigma\cdot (\sigma f)}{2}d\sigma\\
& &+\int_\Sigma(\Delta_\Sigma f)\frac{\sigma\cdot\nabla_\sigma \bar{f}+\nabla_\sigma\cdot (\sigma \bar{f})}{2}d\sigma\\
&=&\sum_{k=1}^{n-1}\int_\Sigma(-\partial_{\sigma_k}^2 \bar{f})(\sigma_k\partial_{\sigma_k} f+\frac{f}{2})d\sigma\\
& &+\int_\Sigma(-\partial_{\sigma_k}^2 f)(\sigma_k\partial_{\sigma_k} \bar{f}+\frac{\bar{f}}{2})d\sigma.
\end{eqnarray*}
Since the computation depends on $\Sigma$, we will only explain here the computation in two particular cases of $\Sigma$ of dimension $2$ (the rectangle and the unitary disc). For other cases, computations are quite similar.

Suppose that $\Sigma=[a,b]\times[c,d]$:
\begin{eqnarray*}
(f,[\Delta_D^\Sigma,iA_D^\sigma]f)&=&\sum_{k=1}^{2}\int_a^b\left(\int_c^d(-\partial_{\sigma_k}^2 \bar{f})(\sigma_k\partial_{\sigma_k} f+\frac{f}{2})d\sigma_2\right)d\sigma_1\\
& &+\int_a^b\left(\int_c^d((-\partial_{\sigma_k}^2 f)(\sigma_k\partial_{\sigma_k} \bar{f}+\frac{\bar{f}}{2})d\sigma_2\right)d\sigma_1.
\end{eqnarray*}
By Fubini Theorem and by integration by part, we have 
 \begin{eqnarray*}
(f,[\Delta_D^\Sigma,iA_D^\sigma]f)&=&\sum_{k=1}^{2}\int_a^b\left(\int_c^d((-2\partial_{\sigma_k}^2 f)\bar{f}d\sigma_2\right)d\sigma_1\\
& &-\int_c^d (b|\partial_{\sigma_1} f|^2(b,\sigma_2)-a|\partial_{\sigma_1} f|^2(a,\sigma_2))d\sigma_2\\
& &-\int_a^b (d|\partial_{\sigma_2} f|^2(\sigma_1,d)-c|\partial_{\sigma_2} f|^2(\sigma_1,c))d\sigma_1.
\end{eqnarray*}
By sum, we obtain for $g\in\cd(\Delta_D)\cap\cd(A_1)$
\begin{eqnarray*}
(g,[\Delta_D,iA_1]g)&=&2(g,\Delta_D g)\\
& &-\int_\R \int_c^d \left(b|\partial_{\sigma_1} g|^2(y,b,\sigma_2)-a|\partial_{\sigma_1} g|^2(y,a,\sigma_2)\right)d\sigma_2 dy\\
& &-\int_\R\int_a^b \left(d|\partial_{\sigma_2} g|^2(y,\sigma_1,d)-c|\partial_{\sigma_2} g|^2(y,\sigma_1,c)\right)d\sigma_1dy.
\end{eqnarray*}

Assume now that $\Sigma=\{(\sigma_1,\sigma_2),\sigma_1^2+\sigma_2^2\leq1\}$ the unitary disc of $\R^2$. Then
\begin{eqnarray*}
(f,[\Delta_D^\Sigma,iA_D^\sigma]f)&=&\sum_{k=1}^{2}\int_\Sigma(-\partial_{\sigma_k}^2 \bar{f})(\sigma_k\partial_{\sigma_k} f+\frac{f}{2})d\sigma\\
& &+\int_\Sigma((-\partial_{\sigma_k}^2 f)(\sigma_k\partial_{\sigma_k} \bar{f}+\frac{\bar{f}}{2})d\sigma.
\end{eqnarray*}
For the term where $k=1$, by Fubini Theorem, we can write
\begin{eqnarray*}
& &\int_\Sigma(-\partial_{\sigma_1}^2 \bar{f})(\sigma_1\partial_{\sigma_1} f+\frac{f}{2})d\sigma\\
&=&\int_{-1}^1\left(\int_{-(1-\sigma_2^2)^{1/2}}^{(1-\sigma_2^2)^{1/2}}(-\partial_{\sigma_1}^2 \bar{f})(\sigma_1\partial_{\sigma_k} f+\frac{f}{2})d\sigma_1
\right)d\sigma_2.
\end{eqnarray*}
Thus, by integration by part, we have
 \begin{eqnarray*}
& &(f,[\Delta_D^\Sigma,iA_D^\sigma]f)\\
&=&\sum_{k=1}^{2}\int_\Sigma(-2\partial_{\sigma_k}^2 f)\bar{f}d\sigma\\
& &-\int_{-1}^1 \biggl((1-\sigma_2^2)^{1/2}|\partial_{\sigma_1} f|^2((1-\sigma_2^2)^{1/2},\sigma_2)\\
& &+(1-\sigma_2^2)^{1/2}|\partial_{\sigma_1} f|^2(-(1-\sigma_2^2)^{1/2},\sigma_2)\biggr)d\sigma_2\\
& &-\int_{-1}^1 \biggl((1-\sigma_1^2)^{1/2}|\partial_{\sigma_2} f|^2(\sigma_1,(1-\sigma_1^2)^{1/2})\\
& &+(1-\sigma_1^2)^{1/2}|\partial_{\sigma_2} f|^2(\sigma_1,-(1-\sigma_1^2)^{1/2})\biggr)d\sigma_1.
\end{eqnarray*}
Therefore, by sum, we obtain for $g\in\cd(\Delta_D)\cap\cd(A_1)$
\begin{eqnarray*}
(g,[\Delta_D,iA_1]g)&=&2(g,\Delta_D g)\\
& &-\int_\R \int_{-1}^1 \biggl((1-\sigma_2^2)^{1/2}|\partial_{\sigma_1} g|^2(y,(1-\sigma_2^2)^{1/2},\sigma_2)\\
& &+(1-\sigma_2^2)^{1/2}|\partial_{\sigma_1} g|^2(y,-(1-\sigma_2^2)^{1/2},\sigma_2)\biggr)d\sigma_2dy\\
& &-\int_\R\int_{-1}^1 \biggl((1-\sigma_1^2)^{1/2}|\partial_{\sigma_2} g|^2(y,\sigma_1,(1-\sigma_1^2)^{1/2})\\
& &+(1-\sigma_1^2)^{1/2}|\partial_{\sigma_2} g|^2(y,\sigma_1,-(1-\sigma_1^2)^{1/2})\biggr)d\sigma_1dy.
\end{eqnarray*}

In the two cases, we can remark that all boundary terms can be seen as an integration on a part of the boundary of $\Omega$ of the function $(y,\sigma_1,\sigma_2)\mapsto \sigma_1|\partial_{\sigma_1} g|^2(y,\sigma_1,\sigma_2)$ or of the function $(y,\sigma_1,\sigma_2)\mapsto \sigma_2|\partial_{\sigma_2} g|^2(y,\sigma_1,\sigma_2)$. For example the term
\[\int_\R \int_{-1}^1 \biggl((1-\sigma_2^2)^{1/2}|\partial_{\sigma_1} g|^2(y,(1-\sigma_2^2)^{1/2},\sigma_2)d\sigma_2dy,\]
which appears in the computation of the commutator when $\Sigma$ is the unitary disc $\R^2$, can be seen as the integral on the set $\{(y,\sigma_1,\sigma_2)\in\R^3\backslash \sigma_1^2+\sigma_2^2=1,\sigma_1\geq0\}$ of the function $(y,\sigma_1,\sigma_2)\mapsto \sigma_1|\partial_{\sigma_1} g|^2(y,\sigma_1,\sigma_2)$, with the parametrisation of this set given by $\{(y,(1-\sigma_2)^{1/2},\sigma_2)\backslash \sigma_2\in[-1,1]\}$.
Thus, we can bounded from above in absolute value all boundary terms by the integral on $\partial\Omega$ of this two functions. Using that the trace of an function is a continuous operator from $\ch^1(\Omega)$ into $L^2(\partial\Omega)$, since $\Sigma$ is bounded by assumptions, we can see that this terms are bounded in norm $\ch^2(\Omega)$. In particular, the commutator is bounded from $\cd(\Delta_D)$ to its dual space and thus $\Delta_D\in C^1(A_1)$. 

Since the regularity $C^1$ is satisfied, we can try to know if the Mourre estimate is true for Dirichlet Laplacian with $A_1$ as conjugate operator. Let
\[H(m)=\begin{cases}
\Delta_y+m E_{\Delta_D^\Sigma}(\{m\})\text{ si }m\in(\lambda_k)_{k\in\N^*};\\
0\text{ sinon}
\end{cases}.\]
Since the application $H(\cdot)$ is equal to zero almost everywhere, we can deduce that $(H(\cdot)+i)^{-1}$ is measurable. Thus, we can write the direct integral
\[\int^\oplus_\R H(m) dm=\Delta_y+\Delta_D^\Sigma=\Delta_D.\]
For an interval $I\subset\R$, this decomposition permits to write the spectral measure of $\Delta_D$ on $I$ with the following form
\[E_{\Delta_D}(I)=\int^\oplus_\R  E_{\Delta_y}(I_m)E_{\Delta_D^\Sigma}(\{m\})dm,\]
with $I_m=\{z\in\R,z+m\in I\}$. Remark that since $\Delta_D^\Sigma$ has a compact resolvent, terms of the previous integral are all equal to zero exepted those for which $m\in (\lambda_k)_{k\in\N}$. In particular, we can write this integral as the following sum
\[E_{\Delta_D}(I)=\sum_{k=1}^\infty  E_{\Delta_y}(I_{\lambda_k})\otimes E_{\Delta_D^\Sigma}(\{\lambda_k\}).\]
Remark that if $I$ is bounded, since $\Delta_y$ is non-negative, $E_{\Delta_y}(I_{\lambda_k})=0$ for $k$ large enough (as soon as $I_{\lambda_k}\subset(-\infty,0)$). In particular, when $I$ is bounded, the previous sum is a finite sum.

If we denote $\psi_k$ eigenvectors of $\Delta_D^\Sigma$, for all $f\in\cd(\Delta_D)$, there is $f_k\in L^2_y$ such that
\begin{equation}\label{eq:decomp f}
f(y,\sigma)=\sum_{k=1}^\infty f_k(y)\psi_k(\sigma).
\end{equation}
Let $I$ a bounded closed interval. We will try to know if the Mourre estimate on $I$ with $A_1$ as conjugate operator is true. Let $f\in\cd(\Delta_D)\cap\cd(A_1)$. We have:
\begin{eqnarray}\label{eq:Mourre est Dirichlet}
(f,E_{\Delta_D}(I)[\Delta_D,\rmi A_1]E_{\Delta_D}(I)f)&=&(f,E_{\Delta_D}(I)[\Delta_y,\rmi A_D^y]E_{\Delta_D}(I)f)\nonumber\\
& &+(f,E_{\Delta_D}(I)[\Delta_D^\Sigma,\rmi A_D^\sigma]E_{\Delta_D}(I)f).
\end{eqnarray}
Using that $[\Delta_y,\rmi A_D^y]=2\Delta_y$ and the decomposition \eqref{eq:decomp f}, we have:
\begin{eqnarray*}
& &(f,E_{\Delta_D}(I)[\Delta_D,\rmi A_1]E_{\Delta_D}(I)f)\\
&=&2\sum_{k,l=1}^\infty\biggl((E_{\Delta_y}(I_{\lambda_k})f_k)\otimes\psi_k,\Delta_y (E_{\Delta_y}(I_{\lambda_l})f_l)\otimes\psi_l\biggr)\\
& &+\sum_{k,l=1}^\infty\biggl((E_{\Delta_y}(I_{\lambda_k})f_k)\otimes\psi_k,[\Delta_D^\Sigma,\rmi A_D^\sigma](E_{\Delta_y}(I_{\lambda_l})f_l)\otimes\psi_l\biggr).
\end{eqnarray*}
Remark again that, since $\Delta_y$ is non-negative, $E_{\Delta_y}(I_{\lambda_k})=0$ for all $k$ large enough and thus sums are finite.

Since $\psi_k$ is an orthonormal family, the first term of the right hand side can be written:
\begin{eqnarray*}
& &2\sum_{k,l=1}^\infty\biggl((E_{\Delta_y}(I_{\lambda_k})f_k)\otimes\psi_k,\Delta_y (E_{\Delta_y}(I_{\lambda_l})f_l)\otimes\psi_l\biggr)\\
&=&2\sum_{k=1}^\infty \int_R(E_{\Delta_y}(I_{\lambda_k})\bar{f}_k)(y)\Delta_y (E_{\Delta_y}(I_{\lambda_k})f_k)(y)dy\\
&=&2\sum_{k=1}^\infty\left(E_{\Delta_y}(I_{\lambda_k}) f_k\otimes\psi_k,\Delta_y E_{\Delta_y}(I_{\lambda_k}) f_k\otimes\psi_k\right)\\
&\geq& 2\sum_{k=1}^\infty\inf(I_{\lambda_k})\left(E_{\Delta_y}(I_{\lambda_k}) f_k\otimes\psi_k, E_{\Delta_y}(I_{\lambda_k}) f_k\otimes\psi_k\right).
\end{eqnarray*}
In particular, if $I$ does not contain any $\lambda_j$, then $0\notin I_{\lambda_k}$ which implies $\inf(I_{\lambda_k})\leq 0$ if and only if $E_{\Delta_y}(I_{\lambda_k})=0$. Taken $a=\min_{k\in\N^*,I_{\lambda_k}\subset(0,+\infty)}\{\inf(I_{\lambda_k})\}>0$,we can show the following
\[E_{\Delta_D}(I)\Delta_y E_{\Delta_D}(I)\geq aE_{\Delta_D}(I).\]

It remains to treat the second part of the right hand side of \eqref{eq:Mourre est Dirichlet}. Since functions $f_k$ do not depend of the variable $\sigma$, we have:
\begin{eqnarray*}
& &\biggl((E_{\Delta_y}(I_{\lambda_k})f_k)\otimes\psi_k,[\Delta_D^\Sigma,\rmi A_D^\sigma](E_{\Delta_y}(I_{\lambda_l})f_l)\otimes\psi_l\biggr)\\
&=&\int_\R (E_{\Delta_y}(I_{\lambda_k})\bar{f}_k)(y)(E_{\Delta_y}(I_{\lambda_l})f_l)(y) dy \cdot\int_\Sigma \bar{\psi}_k(\sigma)[\Delta_D^\Sigma,\rmi A_D^\sigma]\psi_l(\sigma)d\sigma.
\end{eqnarray*}
To begin remark that if there is $j\in\N^*$ such that $I\subset (\lambda_j,\lambda_{j+1})$ and $|I|\leq \inf\limits_{1\leq k\leq j}\lambda_{k+1}-\lambda_k$, then intervals $I_{\lambda_k}$ are separated. This implies that if $k\not=l$, $(E_{\Delta_y}(I_{\lambda_k})f_k)$ et $(E_{\Delta_y}(I_{\lambda_l})f_l)$ are orthogonal. Thus, it remains only  diagonal terms. Since for all $k$, $\psi_k\in\cd(\Delta_D^\Sigma)\subset\cd(A_D^\sigma)$, we obtain:
\begin{eqnarray*}
& &\int_\Sigma \bar{\psi}_k(\sigma)[\Delta_D^\Sigma,\rmi A_D^\sigma]\psi_k(\sigma)d\sigma\\
&=&\int_\Sigma \left(\Delta_D^\Sigma\bar{\psi}_k\right)(\sigma)\rmi \left(A_D^\sigma\psi_k\right)(\sigma)d\sigma-\int_\Sigma \left(A_D^\sigma\bar{\psi}_k\right)(\sigma)\rmi \left(\Delta_D^\Sigma\psi_k\right)(\sigma)d\sigma\\
&=&\lambda_k\int_\Sigma \bar{\psi}_k(\sigma)\rmi \left(A_D^\sigma\psi_k\right)(\sigma)d\sigma-\lambda_k\int_\Sigma \left(A_D^\sigma\bar{\psi}_k\right)(\sigma)\rmi \psi_k(\sigma)d\sigma\\
&=&0.
\end{eqnarray*}
Thus the second term will not neither give some positivity of the commutator nor prevent the positivity. Thus, we have
\begin{eqnarray}\label{eq:commut dirichlet A1}
& &(f,E_{\Delta_D}(I)[\Delta_D,\rmi A_1]E_{\Delta_D}(I)f)\nonumber\\
&=&2\sum_{k,l=1}^\infty\biggl((E_{\Delta_y}(I_{\lambda_k})f_k)\otimes\psi_k,\Delta_y (E_{\Delta_y}(I_{\lambda_l})f_l)\otimes\psi_l\biggr)\\
&\geq&2a\|E_{\Delta_D}(I)f\|^2,\nonumber
\end{eqnarray}
which implies that the Mourre estimate is true on $I$ when $I$ does not contains any $\lambda_j$ and $|I|\leq \inf\limits_{1\leq k\leq n}\lambda_{k+1}-\lambda_k$.

To use the Mourre theorem, we have to proove more regularity of $\Delta_D$ with respect to $A_1$. As for the case where $A$ is a self-adjoint operator, we can define the class of regularity $C^1_u(A)$ if $A$ is maximal symmetric:
\begin{definition}
Let $S$ be a bounded operator and $A$ a maximal symmetric operator such that $S\in C^1(A)$. We say that $S\in C^1_u(A)$ if $\left(Se^{iAt}-e^{iAt}S\right)(it)^{-1}$ has $[S,iA]$ as norm limit.

If $S$ is unbounded, we say that $S\in C^1_u(A)$ if and only if for $z\in\rho(S)$, $(S-z)^{-1}\in C^1_u(A)$
\end{definition}

Now we will prove that the Dirichlet Laplacian is not of class $C^1_u(A_1)$ by using a proof by contradiction. To simplify computations, we will only give details for the case $\Sigma=[-1,1]$ for which eigenvalues and eigenvectors of $\Delta_D^\Sigma$ are well known (see p.266 of \cite{RS4}), but similar proof can be used for other types of $\Sigma$.
In the case $\Sigma=[-1,1]$, eigenvalues are $\lambda_k=\left(\frac{\pi}{2}\right)^2k^2$ with $k\in\N^*$ and associated eigenvectors are
\[\psi_k(\sigma)=\begin{cases}
\cos(k\pi\sigma/2)\text{ si } k\in 2\N+1,\\
\sin(k\pi\sigma/2)\text{ si } k\in 2\N^*.
\end{cases}.\]

By a simple computations, we have for $k,l\in\N^*$
\[
\int_\Sigma \bar{\psi}_k(\sigma)[\Delta_D^\Sigma,\rmi A_D^\sigma]\psi_l(\sigma)d\sigma=\begin{cases}
2\lambda_k\lambda_l(-1)^{\frac{k+l}{2}+1} \text{ si }k-l\in 2\Z^*,\\
0\text{ sinon. }
\end{cases}\]

If we suppose that $\Delta_D\in C^1_u(A_1)$, then, since $\jap{q_y}^{-1}(\Delta_D+1)^{-1}\jap{q_y}^{-1}$ is compact,  $\jap{q_y}^{-1}[(\Delta_D+1)^{-1}, i A_1]\jap{q_y}^{-1}=\jap{q_y}^{-1}(\Delta_D+1)^{-1}[\Delta_D, i A_1](\Delta_D+1)^{-1}\jap{q_y}^{-1}$ is a compact operator, as a limit in norm of compact operators.

Let $g\in\ch^2_1(\R)$, $g\not=0$. Let $h_m(y,\sigma)=g(y)(\psi_{4m}-\psi_{4m+2})(\sigma)$. Since $(\psi_k)$ is an orthonormal family of eigenvectors of $\mathbb{1}\otimes\Delta_\Sigma$, $(h_m)$ is an orthogonal family in $\cg=\ch_{1,y}\cap\cd(\Delta_D)$ with the norm $\|f\|_\cg=\|\jap{q_y}(\Delta_D+1)f\|$, where $\ch_{1,y}$ is the space define by $\|f\|_{\ch_{1,y}}=\|\jap{q_y}f\|_{L^2}$. Thus, we have $[\Delta_D, i A_1]:\cg\rightarrow\cg^*$ compact.  Moreover,
\begin{eqnarray*}
\|h_m\|_\cg^2&=&2\int_\R\jap{y}^2 g''(y)^2dy+2\int_\R\jap{y}^2 g''(y)g(y)dy\biggl(\lambda_{4m}+\lambda_{4m+2}+2\biggr)\\
& &+\int_\R\jap{y}^2 g(y)^2dy\biggl(2+2\lambda_{4m}+2\lambda_{4m+2}+\lambda_{4m}^2+\lambda_{4m+2}^2\biggr).
\end{eqnarray*}
When $m$ goes to infinity, we can observe the equivalence 
\begin{equation}\label{eq: norme hn}
\|h_m\|^2_\cg\simeq\int_\R\jap{y}^2 g(y)^2dy\biggl(\lambda_{4m}^2+\lambda_{4m+2}^2\biggr).
\end{equation}
Let $f_m=h_m/\|h_m\|_\cg$. By definition, $(f_m)$ is an orthonormal family of $\cg$. Thus we have
\begin{eqnarray*}
& &4\int_\R|g'(y)|^2dy\|h_m\|_\cg^{-2}-4\lambda_{4m}\lambda_{4m+2}\|g\|^2_{L^2_y}\|h_m\|_\cg^{-2}\\
&=&(f_m,[\Delta_D,i  A_1]f_m)\\
&\geq&-\|\jap{q_y}^{-1}(\Delta_D+1)^{-1}[\Delta_D,i A_1]f_m\|.
\end{eqnarray*}

Since $(f_m)$ is an orthonormal family in $\cg$, the right hand side goes to $0$ when $m$ goes to infinity. Moreover, since $\lim\limits_{m\rightarrow\infty}\|h_m\|_\cg^{-2}=0$, we deduce that $\lambda_{4m}\lambda_{4m+2}\|g\|^2_{L^2_y}\|h_m\|_\cg^{-2}$ goes to $0$. 

By a simple computation, we can show that 
\[\lim\limits_{m\rightarrow\infty}\frac{\lambda_{4m}\lambda_{4m+2}}{\lambda_{4m}^2+\lambda_{4m+2}^2}=\frac{1}{2}.\]
Thus, using \eqref{eq: norme hn},
\[\lim\limits_{m\rightarrow\infty}\lambda_{4m}\lambda_{4m+2}\|g\|^2_{L^2_y}\|h_n\|_\cg^{-2}=\frac{\|g\|^2_{L^2_y}}{2\|\jap{q_y}g\|^2_{L^2_y}}\not=0.\]

Thus, by contradiction, we deduce that $\Delta_D\notin C^1_u(A_1)$.

Moreover, defining the space $C^{1,1}(A_1)=(C^2_u(A_1),C^0_u(A_1))_{1/2,1}$, we obtain the inclusion $C^{1,1}(A_1)\subset C^1_u(A_1)$. Thus, the Dirichlet Laplacian is not of class $C^{1,1}(A_1)$, a class of regularity which is necessary for Mourre theorem.

Thus Proposition \ref{prop: A_D guide d'onde } is proved.

We can remark that if $I$ contains $\lambda_j$, by taking $f=f_j\otimes\psi_j$, using \eqref{eq:commut dirichlet A1}, we have
\[(f,E_{\Delta_D}(I)[\Delta_D,\rmi A_1]E_{\Delta_D}(I)f)=\|E_{\Delta_y}(I_{\lambda_j})\partial_y f_j\|^2_{L^2(\R)}.\]
since this term is not positive ($0\in I_{\lambda_j}$), we have a problem of threshold, as when we use $A_D^y$ conjugate operator. To solve this problem using the method of the weakly conjugate operator (see \cite{BoGo}), we can try to see if the commutator is non-negative and injective.

To simplify computations, we will only give the details for the case $\Sigma=[-1,1]$. 

Let $g\in \ch^2(\R)$ and $k\in\N$. Let $f=g\otimes\psi_{k+2}-(-1)^k g\otimes \psi_k$. By \eqref{eq:Mourre est Dirichlet}, we have:
\begin{eqnarray}\label{eq: positivite globale}
& &(f,[\Delta_D,\rmi A_1]f)\nonumber\\
&=&2\int_\R \bar{g}(y)\Delta_y g(y)dy-2(-1)^k(g\otimes\psi_{k+2},[\Delta_D^\Sigma,\rmi A_D^\sigma]g\otimes\psi_k)\nonumber\\
&=&2\int_\R  |g'(y)|^2dy-4\lambda_k\lambda_{k+2}\int_\R |g(y)|^2dy.
\end{eqnarray}
Let $h\in\ch^2(\R)$ and let $g_w(y)=w h(w^2y)$ for all $y\in\R$. We can remark that
\[\int_\R |g_w(y)|^2dy=\int_R|h(y)|^2 dy \text{ and }\int_\R  |g_w'(y)|^2dy=w^4\int_\R  |h'(y)|^2dy.\] 
Thus replacing $g$ in \eqref{eq: positivite globale} by $g_w$ and making $w$ go to $0$, we can show that $(f,[\Delta_D,\rmi A_1]f)$ is negative for $w$ small enough. Thus, the commutator is not non-negative (neither injective) which prevent us to use the method of the weakly conjugate operator.

\section{The case of the curved waveguide}\label{s: A_u guide d'onde}

In this section, we will prove a Limiting Absorption Principle for Schr\"odinger operators on curved waveguides. In the following, we will always suppose that $n\geq2$ and we will always be in the context of the article \cite{KT}. The purpose here is to generalize Theorem 3.4 from \cite{KT}, by limiting conditions on derivatives of the curvature of the waveguide. 

\subsection{Geometric preliminaries}

To begin, we will recall some notions concerning geometric properties of waveguides, using notations of \cite{KT}.

Let $p:\R\rightarrow\R^n$ a function of class $C^\infty$. Assume that 
\begin{assumption}\label{a: base}
there is a collection $(e_k)_{k=1,\cdots n}$ of smooth mapping from $\R$ to $\R^n$ such that
\begin{enumerate}[(i)]
\item For all $y\in\R$, $(e_k(y))$ is an orthonormal family;

\item For all $k=1,\cdots,n-1$ and for all $y\in\R$, the $k^{th}$ derivative of $p(y)$ lies in the span of $e_1(y),\cdots, e_k(y)$;

\item $e_1=p'$;

\item For all $y\in\R$, the family $(e_k(y))$ has the positive orientation;

\item For all $k=1,\cdots,n-1$ and for all $y\in\R$, $e_k'(y)$ lies in the span of \\
$e_1(y),\cdots, e_{k+1}(y)$.
\end{enumerate}
\end{assumption}
By the Serret-Frenet formula, we know that there exist a matrix $\ck$ with size $n\times n$ such that 
\[\frac{\partial}{\partial y}(e_k(y))=\ck (e_k(y)).\]
Moreover $\ck$ satisfies:
\[\ck_{i j}=\begin{cases}
\kappa_i\quad\text{ si }j=i+1\\
-\kappa_i\quad\text{ si }i=j+1\\
0\quad\text{ sinon}
\end{cases},\]
where $\kappa_i$ is the $i^{th}$ curvature of $p$.

From $\ck$, we can define the matrix $n\times n$ of rotation $\cR$ which satisfies $\cR_{1 k}=\cR_{k 1}=\delta_{1 k}$ and for all $i,j=2,\cdots,n$, we have
\[\frac{\partial}{\partial s}\cR_{i j}+\sum_{\alpha=1}^n\cR_{i \alpha}\ck_{\alpha j}=0.\]
From the matrix $\cR$,we can define the family $(\tilde{e_k})$ by $(\tilde{e_k})=\cR(e_k)$.

Let $\Sigma$ a bounded open set of $\R^{n-1}$ and $\Omega$ the straight waveguide $\R\times\Sigma$. We can define $\Gamma$ as the image of $\Omega$ by the application 
\[
\cl: \Omega\rightarrow\R^n\quad
(y,\sigma_2,\cdots,\sigma_n)\mapsto p(y)+\sum_{k=2}^n\sigma_k\tilde{e_k}.
\]
Assume that $\Gamma$ does not overlap and that $\kappa_1\in L^\infty$ with $\sup\limits_{\sigma\in\Sigma}|\sigma|\|\kappa_1\|_\infty<1$. Then, $\cl:\Omega\rightarrow\Gamma$ is a diffeomorphism which permits to identify $\Gamma$ with a Riemannian manifold $(\Omega,g)$. Moreover, $g=diag(h^2,1,\cdots,1)$ with 
\[h(y,\sigma)=1+\sum_{k=2}^n\sum_{\alpha=1}^n\sigma_k\cR_{k \alpha}(y)\ck_{\alpha 1}(y)=1-\sum_{k=2}^n\sigma_k\cR_{k 2}(y)\kappa_1(y).\]
We will always assume that $h$ is bounded from below by a positive constant  which implies that the Riemannian metric $g$ is inversible. 

Now it remains to know what operator on $L^2(\Omega,g)$ is associated to the Schr\"odinger operator $\Delta+\tilde{V}$ on $L^2(\Gamma)$ with $\tilde{V}$ a potential on the waveguide $\Gamma$. Remark that, identifying $L^2(\Gamma)$ and $L^2(\Omega,g)$, if we denote $dv$ a volume element of $\Gamma$, this operator is associated to the sesquilinear form $\tilde{Q}$ defined by
\[\tilde{Q}(\phi,\psi)=\int_\Omega g_{1 1}^{-1}\bar{\partial_{y} \phi}\partial_{y}\psi dv+\sum_{k=2}^n\int_\Omega g_{k k}^{-1}\bar{\partial_{\sigma_k} \phi}\partial_{\sigma_k}\psi dv+\int_\Omega \bar{\phi}V\psi dv\]
with appropriate domain ($\ch^1_0$ for Dirichlet boundary conditions, $\ch^1$ for Neumann boundary conditions) and with $V=\cl \tilde{V}\cl$. Since the form $\tilde{Q}$ is densly defined, symmetric and closed on its domain, we can associated to it a unique self-adjoint operator $\tilde{H}$ defined, for $\psi\in\cd(\tilde{H})$, by
\[\tilde{H}\psi=-|g|^{1/2}\left(\partial_{y} |g|^{1/2} g_{1 1}^{-1}\partial_{y}\psi+\sum_{k=2}^n\partial_{\sigma_k} |g|^{1/2} g_{k k}^{-1}\partial_{\sigma_k}\psi\right)+V\psi.\]
To simplify computations, we transform $\tilde{H}$ into a unitary equivalent operator $H$. To do this, we use the unitary transformation $\cu:\psi\mapsto|g|^{1/4}\psi$. By defining $H=\cu\tilde{H}\cu^{-1}$, we obtain
\begin{equation}\label{eq: H guide tordu}
H\psi=-\partial_{y} g_{1 1}^{-1}\partial_{y}\psi-\sum_{k=2}^n\partial_{\sigma_k} g_{k k}^{-1}\partial_{\sigma_k}\psi+(V+W)\psi
\end{equation}
with
\[W=-\frac{5}{4}\frac{(\partial_y h)^2}{h^4}+\frac{1}{2}\frac{\partial_y^2 h}{h^3}-\frac{1}{4}\frac{\sum_{k=2}^n(\partial_{\sigma_k} h)^2}{h^2}+\frac{1}{2}\frac{\sum_{k=2}^n\partial_{\sigma_k}^2 h}{h}.\]
Remark that with our choice of $h$, for all $k=2,\cdots,n$, $\partial_{\sigma_k}^2 h=0$ and $g_{k k}=1$.

\subsection{A Limiting Absorption Principle far from threshods}
 
Now, we will prove a Limiting Absorption Principle far from thesholds. To do this, we will use Mourre theorem and for this reason, we have to find a conjugate operator. As we saw previously (c.f. Section \ref{s: A_D toute direction}), it seems necssary to take for conjugate operator an operator only in the unbounded direction of the waveguide.

As for the Euclidian space $\R^n$, a natural conjugate operator to use is the generator of dilations (see \cite{CFKS,Mo1,Mo2,KT}). To apply Mourre theorem, it is sufficient to assume the following: :
\begin{assumption}[Assumption 3.3 de \cite{KT}]
Uniformly in $\sigma\in\Sigma$,
\begin{enumerate}
\item $h(y,\sigma)\rightarrow 1$ when $|y|\rightarrow\infty$;

\item $\partial^2_yh(y,\sigma),\sum_{k=2}^n(\partial_{\sigma_k}h(y,\sigma))^2\rightarrow 0$ when $|y|\rightarrow\infty$;

\item there is $\theta\in(0,1]$ such that
\[\partial_yh(y,\sigma),\partial^3_yh(y,\sigma),\sum_{k=2}^n\partial_y(\partial_{\sigma_k}h)^2(y,\sigma)=O(|y|^{-1-\theta}).\]
\end{enumerate}
\end{assumption}
Under these assumptions and under similar assumptions for the decay of $V$ than in theorem \ref{th: KT}, spectral results of theorem \ref{th: KT} stay true for $H$ a Schr\"odinger operator on $L^2(\Gamma)$. Remark that the fact that $\sigma_{ess}(H)=[\nu,\infty)$, with $\nu$ the first eigenvalue of the Dirichlet Laplacian on $L^2(\Sigma)$ does not depend of the fact that $\partial_y^2 h(y,\sigma)\rightarrow 0$ when $|y|\rightarrow\infty$. In fact, if we remove this assumption, using that 
\[\frac{\partial_y^2 h}{h^3}=\partial_y(\frac{\partial_y h}{h^3})+2\frac{(\partial_y h)^2}{h^4},\]
and writing that $\partial_y(\frac{\partial_y h}{h^3})=[ip_y,\frac{\partial_y h}{h^3}]$, we can show that $W:\ch^1\rightarrow\ch^{-1}$ is compact which does not change the essential spectrum of $H_0=-\partial_{y} g_{1 1}^{-1}\partial_{y}-\sum_{k=2}^n\partial_{\sigma_k} g_{k k}^{-1}\partial_{\sigma_k}$. 

Since derivatives of $h$ can be written as a function of  coefficients of the matrix $\cR$ and of derivatives of coefficients of the matrix $\ck$, all assumptions on derivatives of $h$ can be written as assumptions on curvatures:

\begin{assumption}[Assumption 3.4 de \cite{KT}]
For all $\alpha\in \{2,\cdots,n\}$,
\begin{enumerate}
\item $\ck^1_\alpha(y),\partial_y^2\ck^1_\alpha (y)\rightarrow 0$ when $|y|\rightarrow\infty$;

\item for all $\beta\in\{2,\cdots,n\}$,$\ck_\alpha^\beta,\partial_y\ck^2_\alpha\in L^\infty(\R)$;

\item there is $\theta\in(0,1]$ such that
\begin{multline*}
\partial_y\ck^1_\alpha(y),\partial_y^3\ck^1_\alpha(y),\ck^2_\alpha(y),\partial_y^2\ck^2_\alpha(y),\sum_{k=2}^n\ck^k_\alpha(y)\partial_y\ck^2_k(y),\sum_{k=2}^n\partial_y\ck^k_\alpha(y)\ck^2_k(y)\\
=O(|y|^{-1-\theta}).
\end{multline*}
\end{enumerate}
\end{assumption}
We can see that using the generator of dilations as conjugate operator impose that curvatures have to be regular functions and that some of their derivatives have some decay at infinity.

To generalize this result, we can choose an other conjugate operator which avoid us to take derivative of the potential $W$. Let $\lambda:\R\rightarrow\R$ a $C^\infty$ function, positive, bounded, with all derivatives bounded such that for all $\alpha\in\N$, $y\rightarrow y\partial_y^\alpha \lambda(y)$ is bounded . Let $A_u=\frac{1}{2}(y p_y\lambda(p_y)+\lambda(p_y)p_y y)$ with $p_y=-i\partial_y$. In \cite[Theorem 4.2.3]{ABG}, we can see that $A_u$ is essentially self-adjoint with domain $C^\infty_c(\R)$. Let $H_1=H-V$. We will compute the first commutator define as a form on $\cd(H_1)\cap\cd(A_u)$:
\begin{eqnarray}\label{eq: commut H guide d'onde}
[H_1,iA_u]&=&[p_y,iA_u]h^{-2}p_y+p_y[h^{-2},iA_u]p_y+p_yh^{-2}[p_y,iA_u]\nonumber\\
& &+[W,iA_u]\nonumber\\
&=&p_y\lambda(p_y)h^{-2}p_y+p_yh^{-2}p_y\lambda(p_y)+p_y[h^{-2},iA_u]p_y\nonumber\\
& &+[W,iA_u]\nonumber\\
&=& 2p_y\lambda(p_y)^{1/2}h^{-2}\lambda(p_y)^{1/2}p_y\nonumber\\
& &+p_y[\lambda(p_y)^{1/2},[\lambda(p_y)^{1/2},h^{-2}]]p_y+p_y[h^{-2},iA_u]p_y\nonumber\\
& &+[W,iA_u].
\end{eqnarray}
We can remark that if we suppose that there is a constant $a>0$ such that $h^{-2}(y,\sigma)\geq a$, for all $(y,\sigma)\in\Omega$ then 
\[2p_y\lambda(p_y)^{1/2}h^{-2}\lambda(p_y)^{1/2}p_y\geq0.\]
In particular, far from thresholds (far from $\ct$), this term will give us positivity, necessary to the obtention of Mourre estimate.
 
 Now we will give some sufficient assumptions to have the good regularity
\begin{assumption}\label{a: hypothese sur h}
Let $h$ such that
\begin{enumerate}
\item $h^{-2}$ is bounded.

\item $h(y,\sigma)\rightarrow 1$ when $|y|\rightarrow\infty$;

\item ther is $b>0$ tel que $h^{-2}(y,\sigma)\geq b$ for all $(y,\sigma)\in\Omega$.

\item there is $\theta>0$ such that, iniformly in $\sigma\in\Sigma$, 
\[\sum_{k=2}^n(\partial_k h)^2(y,\sigma)=O(|y|^{-(1+\theta)}) \text{ and } \partial_y h(y,\sigma)=O(|y|^{-(1+\theta)}).\]
\end{enumerate}
\end{assumption}
 Under these assumptions, we have the following result:

\begin{theorem}\label{th: LAP guide}
Let $\Gamma$ a waveguide as it was definite previously. Suppose assumptions \ref{a: base} and \ref{a: hypothese sur h}. Assume moreover, that for all $\alpha\in\N$,  $y\mapsto\jap{y}^{1+\alpha}\partial_y^\alpha\lambda(y)$ is bounded. Let $V$ a potential compact from $\ch_y^1$ to $(\ch^1_y)^*$ of class $C^{1,1}(A_u,\ch_y^1,\ch_y^{-1})$. Then spectral results of  Theorem \ref{th: KT} stay true for $H=\Delta+\tilde{V}$ with Dirichlet boundary conditions with $\tilde{V}=\cl^{-1} V\cl^{-1}$.
\end{theorem}
We can remark that if we assume that $\partial_y^2 h(y,\sigma)$ goes to $0$ when $|y|\rightarrow\infty$, uniformly in $\sigma\in\Sigma$, then, if we suppose that $V$ is $\Delta$-compact and of class $C^{1,1}(A_u,\ch^2,\ch^{-2})$, Theorem \ref{th: LAP guide} stay true.

\begin{proof}[Theorem \ref{th: LAP guide}]
Let $\ch^1_y$ be the domain of $\jap{p_y}$. Remark that the form domain of the Dirichlet Laplacian $\cq(\Delta_D)$ is a subset of $\ch^1_y$.
If $h^{-2}$ is bounded, we can remark that
 \[\jap{p_y}^{-1}p_y\lambda(p_y)^{1/2}h^{-2}\lambda(p_y)^{1/2}p_y\jap{p_y}^{-1}\] and 
\[\jap{p_y}^{-1}p_y[\lambda(p_y)^{1/2},[\lambda(p_y)^{1/2},h^{-2}]]p_y\jap{p_y}^{-1}\] are bounded. Moreover, if $h$ satisfies assumptions \ref{a: hypothese sur h}, by writing 
\[\frac{\partial_y^2 h}{h^3}=\partial_y(\frac{\partial_y h}{h^3})+2\frac{(\partial_y h)^2}{h^4},\]
$\jap{p_y}^{-1}W\jap{p_y}^{-1}$ is compact. With a similar proof, we can show that $\jap{p_y}^{-1}\jap{q}^{1+\theta}W\jap{p_y}^{-1}$ is bounded which implies that $\jap{y}^\theta\jap{p_y}^{-1}[W,iA_u]\jap{p_y}^{-1}$ is bounded, since $p_y\lambda(p_y)$ is bounded. Since 
\[[h^{-2},iA_u]=y[h^{-2},i p_y\lambda(p_y)]+\frac{1}{2}[h^{-2},\lambda(p_y)+p_y\lambda'(p_y)],\]
by the Helffer-Sjostrand formula, we can see that $\jap{y}^\theta[h^{-2},iA_u]$ is bounded. Moreover, using that
\[[\lambda(p_y)^{1/2},[\lambda(p_y)^{1/2},h^{-2}]]=\lambda(p_y)^{1/2}[\lambda(p_y)^{1/2},h^{-2}]-[\lambda(p_y)^{1/2},h^{-2}]\lambda(p_y)^{1/2}\]
 and by the Helffer-Sjostrand formula, we deduce that $\jap{y}[\lambda(p_y)^{1/2},[\lambda(p_y)^{1/2},h^{-2}]]$ is bounded.

By a commutator computation, we have
\begin{eqnarray*}
& &[p_y\lambda(p_y)^{1/2}h^{-2}\lambda(p_y)^{1/2}p_y,i A_u]\\
&=&[p_y\lambda(p_y)^{1/2},i A_u]h^{-2}\lambda(p_y)^{1/2}p_y+p_y\lambda(p_y)^{1/2}[h^{-2},i A_u]p_y\lambda(p_y)^{1/2}\\
& &+p_y\lambda(p_y)^{1/2}h^{-2}[\lambda(p_y)^{1/2}p_y,i A_u]\\
&=&\left(\lambda(p_y)^{1/2}+\frac{1}{2}p_y\partial_y\lambda(p_y)\lambda(p_y)^{-1/2}\right)p_y\lambda(p_y)h^{-2}\lambda(p_y)^{1/2}p_y\\
& &+p_y\lambda(p_y)^{1/2}[h^{-2},i A_u]p_y\lambda(p_y)^{1/2}\\
& &+p_y\lambda(p_y)^{1/2}h^{-2}\left(\lambda(p_y)^{1/2}+\frac{1}{2}p_y\partial_y\lambda(p_y)\lambda(p_y)^{-1/2}\right)p_y\lambda(p_y).
\end{eqnarray*} 
In this way, we can prove that the commutator is bounded from $\ch^1$ to $\ch^{-1}$ which implies that the first term of the right hand side of  \eqref{eq: commut H guide d'onde} is of class \\
$C^1(A_u,\ch_y^1,\ch_y^{-1})\subset C^{0,1}(A_u,\ch_y^1,\ch_y^{-1})$. In particular, this implies that $H_1=H-V$ is of class $C^{1,1}(A_u,\ch_y^1,\ch_y^{-1})$. Since $V$ is a compact potential from $\ch^1$ to $\ch^{-1}$ and since $V$ is of class $C^{1,1}(A_u,\ch_y^1,\ch_y^{-1})$, we deduce by sum that $H$ is of class $C^{1,1}(A_u,\ch_y^1,\ch_y^{-1})$

It remains to prove that the Mourre estimate is satisfied for $H_1$ with $A_u$ as conjugate operator on all compact interval of $(\nu,+\infty)\backslash\ct$ where $\nu=\inf\ct$. If we denote $H'=p_yh^{-2}p_y+\Delta_D^\Sigma$, since $H$ and $H'$ are of class $C^{1,1}(A_u)$, we can remark that it is sufficient to prove that the Mourre estimate is true near all points of $(\nu,+\infty)\backslash\ct$ for $H'$ with $A_u$ as conjugate operator to obtain a Mourre estimate near all points of $(\nu,+\infty)\backslash\ct$ for $H$ with $A_u$ as conjugate operator (see Theorem 7.2.9 of \cite{ABG}). Since $A_u$ can be written $A_1\otimes\mathbb{1}_\Sigma+\mathbb{1}_\R\otimes A_2$ with $A_2=0$, using Theorem 2.9 of \cite{KT}, we can see that it is sufficient to prove that the Mourre estimate is true near all points of  $\R^{+*}$ for $H_0=p_yh^{-2}p_y$ with $A_u$ as conjugate operator to obtain a Mourre estimate near all points of $(\nu,+\infty)\backslash\ct$ for $H$ with $A_u$ as conjugate operator, since $\Delta_D^\Sigma$ has a pure point spectrum.

Let $\lambda\in (0,+\infty)$ and $\phi\in C^\infty_c(\R)$ such that $\phi(\lambda)\not=0$. We will show that there exists $a>0$ and $K$ compact such that
\[\phi(H_0)[H_0,i A_u]\phi(H_0)\geq a\phi(H_0)^2+K.\]
In the following, to simpify notations, we denote $K_k$, $k\in\N$ some compact operators. By \eqref{eq: commut H guide d'onde}, we have
\begin{eqnarray*}
[H_0,iA_u]&=& 2p_y\lambda(p_y)^{1/2}h^{-2}\lambda(p_y)^{1/2}p_y\\
& &+p_y[\lambda(p_y)^{1/2},[\lambda(p_y)^{1/2},h^{-2}]]p_y+p_y[h^{-2},iA_u]p_y.
\end{eqnarray*}
Remark that the two last terms of the right hand side are compact from  $\ch^1_y$ to $\ch^{-1}_y$ which implies that we have
\begin{equation}\label{eq: est Mourre guide tordu}
\phi(H_0)[H_0,i A_u]\phi(H_0)=2\phi(H_0)p_y\lambda(p_y)^{1/2}h^{-2}\lambda(p_y)^{1/2}p_y\phi(H_0)+K_1.
\end{equation}
By a simple computation, we have
\[(H_0+i)^{-1}-(p_y^2+i)^{-1}=(H_0+i)^{-1}p_y(1-h^{-2})p_y(p_y^2+i)^{-1}.\]
Since $h(y,\sigma)\rightarrow 1$ when $|y|\rightarrow\infty$ uniformly in $\sigma\in\Sigma$, this implies that $(H_0+i)^{-1}-(p_y^2+i)^{-1}$ is compact. By Lemma 7.2.8 of \cite{ABG}, we deduce that $\phi(H_0)-\phi(p_y^2)$ is compact. Using \eqref{eq: est Mourre guide tordu}, we thus have
\[\phi(H_0)[H_0,i A_u]\phi(H_0)=2\phi(p_y^2)p_y\lambda(p_y)^{1/2}h^{-2}\lambda(p_y)^{1/2}p_y\phi(p_y^2)+K_2.\]
Let $\epsilon>0$. By choosing $\phi$ such that $\phi(\alpha)=0$ for all $\alpha<\epsilon$, since $h$ est bornée etis bounded and $\lambda$ is positive, we deuce that there is $a>0$ such that 
\[\phi(H_0)[H_0,i A_u]\phi(H_0)\geq a\phi(p_y^2)^2+K_2=a\phi(H_0)^2+K_3.\]
Therefore, the Mourre estimate is true for $H_0$ with $A_u$ as conjugate operator near all points of $\R^{+*}$ which implies that the Mourre estimate is true for $H$ with $A_u$ as conjugate operator near all points of $(\nu,+\infty)\backslash\ct$.

Theorem \ref{th: LAP guide} is then a consequence of Mourre Theorem.
 \qed
\end{proof}

We can remark again that, since the conjugate operator does not depends on the bounded direction of the waveguide, boundary conditions does not appear. In particular, if $\Gamma$ is a waveguide with a boundary with enough regularity (at least $C^1$), we can define at all points of $\partial\Gamma$ a tangent space and thus a normal derivative at the boundary. Since $\cl$ is a diffeomorphism, it sends tangent spaces of $\Gamma$ into tangent spaces of $\Omega$. Theredore, Neumann/Robin boundary conditions on $\Gamma$ are transform into Neumann/Robin boundary conditions on $\Omega$. In particular, with a similar proof, Theorem \ref{th: LAP guide} stays true if we replace Dirichlet boundary conditions by Neumann or Robin boundary conditions.

Using the particular form of $h$, we can translate assumptions \ref{a: hypothese sur h} into assumptions on curvatures $\kappa_k$:

\begin{assumption}\label{a: courbure A_u}
assume that $\kappa_1\in L^\infty$ with $\sup\limits_{\sigma\in\Sigma}|\sigma|\|\kappa_1\|_\infty<1$. Assume moreover that there is $\theta>0$ such that
\begin{enumerate}
\item $\lim\limits_{|y|\rightarrow\infty}\kappa_1(y)=0$;

\item $\kappa_1'(y)=O(|y|^{-(1+\theta)})$;

\item $\kappa_2(y)\kappa_1(y)=O(|y|^{-(1+\theta)})$.
\end{enumerate}
\end{assumption}
We can remark that if $\kappa_1\in L^\infty$ and if $\sup\limits_{\sigma\in\Sigma}|\sigma|\|\kappa_1\|_\infty<1$ then $h^{-2}$ is bounded and there is $b>0$ such that $h^{-2}(y,\sigma)\geq b$ for all $(y,\sigma)\in\Omega$.
Denoting that $\sum_{k=2}^n(\partial_{\sigma_k} h)^2(y,\sigma)=\kappa_1^2$, we can remark that these assumptions implies assumptions \ref{a: hypothese sur h}. Thus we can write Theorem \ref{th: LAP guide} with assumptions on curvatures:
\begin{theorem}
Let $\Gamma$ a waveguide as it was defined previously. Suppose that assumptions \ref{a: base} and \ref{a: courbure A_u} are satisfied. Let $V$ a compact potential from $\ch_y^1$ to $(\ch^1_y)^*$ of class $C^{1,1}(A_u,\ch_y^1,\ch_y^{-1})$. Then spectral results of Theorem \ref{th: KT} are true for $H=\Delta+V$ with Dirichlet boundary conditions.
\end{theorem}

Remark that the use of an operator $A_u$ at the place of the generator of dilations, usually used, permits two improvements. In a first time, we can avoid to impose conditions on $(\kappa_k)_{k=3,\cdots,n}$ and on  derivatives of $\kappa_2$. Moreover, we can remark that if $\kappa_1$ goes quickly to $0$ at infinity, $\kappa_2$ is not necesseraly bounded. Moreover, no conditions on derivatives of order higher or equal to 2 of $\kappa_1$ are impose, we can take as curvatures $\kappa_k$ functions with high oscillations. The second remark we can do is that, as in \cite{Ma1}, the use of $A_u$ permits to treat a larger class of potential, without assuming for example that the potential is a regular function or is $\Delta$-compact.

\subsection{A Limiting Absorption Principle near threshods}
In this section, we will use the version of the Mourre theory called the method of the weakly conjugate operator, to obtain a Limiting Absorption Principle near thresholds. For more details about this method, we refer to papers \cite{BM,Ri,BoGo,Ma3}.

Let $F:\R\rightarrow\R$ of class $C^\infty$ with all derivatives bounded. Let $A=A_F\otimes 1_\Sigma$ where $A_F$ is define by
\[A_F=\frac{1}{2}\left(p_y F(y)+F(y)p_y\right)\]
with $p_y=-i\partial_y$.
In this case, we have on $\cd(H)\cap\cd(A_F)$:
\begin{eqnarray*}
[H,iA]&=&[p_{y} g_{1 1}^{-1}p_{y},iA]+[V+W,iA]\\
&=&[p_{y} ,iA]g_{1 1}^{-1}p_{y}+p_{y}[g_{1 1}^{-1},iA]p_{y}+p_{y} g_{1 1}^{-1}[p_{y},iA]+[V+W,iA].
\end{eqnarray*}

By a commutator computation on the form domain of $\Delta_D$, we can see that
\begin{eqnarray*}
[p_y,iA]&=&\frac{1}{2}\left(p_y[p_y,iF(y)]+[p_y,iF(y)]p_y\right)\\
&=&\frac{1}{2}\left(p_yF'(y)+F'(y)p_y\right)\\
&=&p_yF'(y)+\frac{i}{2}F''(y)\\
&=&F'(y)p_y-\frac{i}{2}F''(y).
\end{eqnarray*}
This implies that
\begin{eqnarray*}
[H,iA]&=&2p_yF'(y)g_{1 1}^{-1}p_{y}+\frac{i}{2}F''(y)g_{1 1}^{-1}p_{y}-\frac{i}{2}p_{y} g_{1 1}^{-1}F''(y)\\
& &+p_{y}[g_{1 1}^{-1},iA]p_{y}+[V+W,iA]\\
&=&2p_yF'(y)g_{1 1}^{-1}p_{y}+\frac{1}{2}[iF''(y)g_{1 1}^{-1},p_{y}]+p_{y}[g_{1 1}^{-1},iA]p_{y}+[V+W,iA]\\
&=&2p_yF'(y)g_{1 1}^{-1}p_{y}-\frac{1}{2}F'''(y)g_{1 1}^{-1}-\frac{1}{2}F''(y)[p_y,ig_{1 1}^{-1}]\\
& &+p_{y}[g_{1 1}^{-1},iA]p_{y}+[V+W,iA].
\end{eqnarray*}
Using that $[G(y),iA]=-F(y)G'(y)$ for all functions $G:\R\rightarrow\R$ and that $g_{1 1}=h^2$, we have 
\begin{eqnarray*}
[H,iA]&=&2p_yF'(y)h^{-2}(y,\sigma)p_{y}-\frac{1}{2}F'''(y)h^{-2}(y,\sigma)+F''(y)\partial_y h(y,\sigma)h^{-3}(y,\sigma)\\
& &+2p_{y}F(y)\partial_y h(y,\sigma)h^{-3}(y,\sigma)p_{y}-F(y)\partial_y(V+W)\\
&=&2p_yh^{-2}\left(F'(y)+F(y)\partial_y h h^{-1}\right)p_{y}\\
& &-\frac{1}{2}F'''(y)h^{-2}+F''(y)\partial_y h h^{-3}-F(y)\partial_y(V+W).
\end{eqnarray*}

To apply the method of the weakly conjugate operator, the commutator has to be non-negative. Thus, assume the following
\begin{assumption}\label{H:guide d'onde A_F}
Assume that for all $(y,\sigma)\in\Omega$,
\begin{enumerate}
\item $F'(y)+F(y)\partial_y h(y,\sigma) h^{-1}(y,\sigma)>0$ is bounded;

\item $(y,\sigma)\mapsto F(y)\partial_y(V+W)(y,\sigma)$ is bounded;
\item $-\frac{1}{2}F'''(y)h^{-2}(y,\sigma)+F''(y)\partial_y h(y,\sigma) h^{-3}(y,\sigma)-F(y)\partial_y(V+W)(y,\sigma)\geq0$ is bounded.
\end{enumerate}
\end{assumption}
Under these assumptions, $S=[H,iA]$ is non-negative and injective. Moreover, $\exp(it A)$ leaves invariant $\cg=\cd(S^{1/2})=\ch^1_y$ (see \cite[Proposition 4.2.4]{ABG}). It remains to prove that $S\in C^1(A_F,\cs,\cs^*)$ where $\cs$ is the completion of $\cg$ for the norm $\|\cdot\|_\cs=\|S^{1/2}\cdot\|$. By a commutator computations, we obtain on $\cd(S)\cap\cd(A_F)$
\begin{eqnarray}\label{eq: commut S guide}
[S,iA_F]&=&4p_y\left(F'(y)h^{-2}(F'(y)+F(y)\partial_y h h^{-1})\right)p_y\nonumber\\
& &-2p_yF(y)\partial_y\left(h^{-2}(F'(y)+F(y)\partial_y h h^{-1})\right)p_y\nonumber\\
& &-\partial_y\left(F''(y)h^{-2}(F'(y)+F(y)\partial_y h h^{-1})\right)\nonumber\\
& &+F(y)\partial_y\left(\frac{1}{2}F'''(y)h^{-2}-F''(y)\partial_y h h^{-3}\right)\nonumber\\
& &+(F(y)\partial_y)^2(V+W).
\end{eqnarray}
We can remark that if $F'$ is bounded, then, the first term of the right hand side is bounded from $\cs$ to $\cs^*$.

By a simple computation, we have:
\begin{eqnarray}\label{eq: commut S guide 2}
& &F(y)\partial_y\left(h^{-2}(F'(y)+F(y)\partial_y h h^{-1})\right)\nonumber\\
&=&-2F(y)\partial_y h h^{-1}\left(h^{-2}(F'(y)+F(y)\partial_y h h^{-1})\right)\nonumber\\
& &+F(y)h^{-2}\left(F''(y)+F'(y)\partial_y h h^{-1}+F(y)\partial_y^2 h h^{-1}-F(y)(\partial_y h)^2 h^{-2}\right).
\end{eqnarray}
Since $F'$ is bounded, we can see that assumptions \ref{H:guide d'onde A_F} imply that $F(y)\partial_y h h^{-1}$ is bounded. In particular, the first term of the right hand side is bounded from above by \\
$\left(h^{-2}(F'(y)+F(y)\partial_y h h^{-1})\right)$. Thus, to get conditions on the second order commutator, it suffices to assume the following:
\begin{assumption}\label{H:guide d'onde 2}
For all $(y,\sigma)\in\Omega$, let 
\[G(y,\sigma)=F'(y)+F(y)\partial_y h(y,\sigma) h^{-1}(y,\sigma)\]
and
\[W_1(y,\sigma)=-\frac{1}{2}F'''(y)h^{-2}(y,\sigma)+F''(y)\partial_y h(y,\sigma) h^{-3}(y,\sigma)-F(y)\partial_y(V+W)(y,\sigma).\]
Assume that for all $(y,\sigma)\in\Omega$, 
\begin{enumerate}
\item there is $C_1>0$ such that
\[
\left|F(y)(F''(y)+F'(y)\partial_y h h^{-1}+F(y)\partial_y^2 h h^{-1}-F(y)(\partial_y h)^2 h^{-2})\right|\leq C_1G(y,\sigma);
\]

\item there is $C_2>0$ such that
\begin{eqnarray*}
& &\biggl|-\partial_y\left(F''(y)h^{-2}(y,\sigma)(F'(y)+F(y)\partial_y h(y,\sigma) h^{-1}(y,\sigma))\right)\\
& &+F(y)\partial_y\left(\frac{1}{2}F'''(y)h^{-2}(y,\sigma)-F''(y)\partial_y h(y,\sigma) h^{-3}(y,\sigma)\right)\\
& &+(F(y)\partial_y)^2(V+W)(y,\sigma)\biggr|\\
&\leq&C_2 W_1(y,\sigma).
\end{eqnarray*}
\end{enumerate}
\end{assumption}
We can remark that $S=2p_y h^{-2}Gp_y+W_1$. Thus, to show that $[S,\rmi A]$ is bounded from $\cs$ to $\cs^*$, it suffices to show that $[S,\rmi A]$ can be bounded from above in the form sense by \\
$Cp_y h^{-2}Gp_y+C' W_1$ with $C,C'>0$ two constants. As it was said previously, since $F'$ is bounded,the first term of the right hand side of \eqref{eq: commut S guide} can be bounded from above in the form sense by $4\|F'\|_{L^\infty}p_y h^{-2}Gp_y$. Moreover, using \eqref{eq: commut S guide 2}, we can remark that, under assumptions \ref{H:guide d'onde 2}, the function $y\mapsto F(y)\partial_y\left(h^{-2}(y)(F'(y)+F(y)\partial_y h(y) h^{-1}(y))\right)$ can be bounded from above by $(2+C_1)Gh^{-2}$. Moreover, by assumptions, the sum of the three last terms of the right hand side of \eqref{eq: commut S guide}  can be bounded from above in the form sense by $C_2 W_1$.
Thus, we can show that $[S,\rmi A]$ is bounded from $\cs$ to $\cs^*$ which implies that $S$ is of class $C^1(A,\cs,\cs^*)$. As it was said previously, this regularity of the operator $S$ implies that $\exp(itA)$ leaves invariant the space $\cs$. Thus $S$ is a good candidate to use the method of the weakly conjugate operator. Moreover, since $S=[H,\rmi A]$, this regularity implies also that $H$ is of class $C^2(A,\cs,\cs^*)$. Therefore, we obtain the following:
\begin{theorem}
Assume that $V\in L^1_{loc}(\R^n,\R)$ is a potential $\Delta$-bounded with bound small than 1. Let $F:\R\rightarrow\R$ of class $C^\infty$ with all derivatives bounded such that assumptions \ref{H:guide d'onde A_F} and \ref{H:guide d'onde 2} are satisfied. Then, there is $c>0$ such that
\[\left|(f,(H-\lambda+\rmi \eta)^{-1}f)\right|\leq \|S^{-1/2}f\|^2+\|S^{-1/2}Af\|^2,\]
with $S=[H,\rmi A]$.

Moreover, $H$ does not have real eigenvalues.
\end{theorem} 

\medskip

{\bf Acknowledgements.} This article is a translation of the last chapter of my PhD thesis (see \cite{MaThese}).
I thank my doctoral supervisor, Thierry Jecko, for fruitful discussion and comments. 

\bibliographystyle{alpha}
\bibliography{bibliographie}

\end{document}